\documentclass[12pt,oneside]{article}
\RequirePackage{epsfig}
\RequirePackage{amsmath}
\RequirePackage{amsfonts,amssymb}
\RequirePackage{amsthm}
\newcommand{\nwt}{\newtheorem}
\nwt{theo}{Theorem}
\nwt{lem}{Lemma}
\nwt{prop}{Proposition}
\nwt{defn}{Definition}
\newcommand{\nwc}{\newcommand}
\nwc{\noi}{\noindent}
\nwc{\non}{\nonumber}
\nwc{\p}{\partial}
\nwc{\half}{\frac{1}{2}}
\nwc{\ca}[1]{{\cal #1}}

\nwc{\T}{\textrm}
\nwc{\wh}{\widehat}
\nwc{\wt}{\widetilde}
\nwc{\bfm}[1]{\mbox{\boldmath $#1$}}
\nwc{\vmm}[1]{\vspace{#1 mm}}
\nwc{\hmm}[1]{\hspace{#1 mm}}
\nwc{\Natural}{\mathbb{N}}

\newcommand{\jumptwo}[1]{[\![#1]\!]}
\nwc{\Ia}{I_\alpha}
\nwc{\Ib}{I_\beta}
\nwc{\cargo}{C}
\nwc{\Rt}{R_t}
\nwc{\Rd}{R_d}
\nwc{\R}{R}
\nwc{\Rtb}{C_{tb}}
\nwc{\Rta}{C_{ta}}
\nwc{\Iab}{C_{ab}}
\nwc{\Rg}{R_g}
\nwc{\Tr}{T}
\nwc{\Iabc}{C_{abc}}
\nwc{\Rc}{C_{1}}
\nwc{\Iac}{C_{a c}}
\nwc{\Kct}{K_c^t}
\nwc{\Kmt}{K_m^t}
\nwc{\Kcd}{K_c^d}
\nwc{\Kmd}{K_m^d}
\nwc{\normale}{{\mathbf{n}}}
\newcommand{\mbf}[1]{\mbox{\boldmath$\rm{#1}$}}
\newcommand{\bu}{{\mbf{u}}}
\newcommand{\bv}{{\mbf{v}}}
\newcommand{\Gatr}{{\Gamma_{\rm tr}}}
\newcommand{\honen}{\mathcal{H}^1}
\newcommand{\be}{\begin{equation}}
\newcommand{\ee}{\end{equation}}

\newcommand{\el}{ \kappa \in \mathcal{T}}
\newcommand{\dint}{\text{\rm int}}
\DeclareMathOperator{\diag}{diag}
\newcommand{\mean}[1]{ \{#1\} }
\newcommand{\ha}{\frac{1}{2}}
\newcommand{\jumpthree}[1]{[\![\![#1]\!]\!]}
\newcommand{\tth}{\mbf{\tt h}}
\DeclareMathOperator{\diam}{diam}
\newcommand{\ltwoin}[2]{\langle{#1},{#2}\rangle}
\newcommand{\texte}[1]{\quad \text{{#1}}\quad}
\newcommand{\su}{\sum_{\kappa\in\mathcal{T}}}
\newcommand{\Gint}{\Gamma_{\dint}}
\newcommand{\bea}{\begin{eqnarray}}
\newcommand{\eea}{\end{eqnarray}}
\newcommand{\beaa}{\begin{eqnarray*}}
\newcommand{\eeaa}{\end{eqnarray*}}
\begin{document}
\thispagestyle{empty}

\begin{center}
\begin{flushleft}
\begin{minipage}[l]{11.5cm}
\textsc{\large Universit\`{a} di Milano--Bicocca\\
Quaderni di Matematica}
\end{minipage}
\hfill
\begin{minipage}[r]{2.5cm}
\includegraphics[height=2cm]{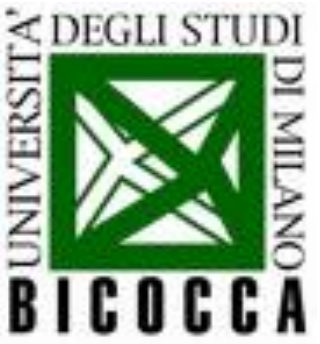}
\end{minipage}
\end{flushleft}

\vspace{2cm}

\textsc{\large Modeling and simulation of Ran-mediated nuclear import}

\vspace{.8cm}

\textsc{Andrea Cangiani}

\vfill

\textsc{Quaderno n. 10/2009} 

\vspace{1cm}

\includegraphics[width=8cm]{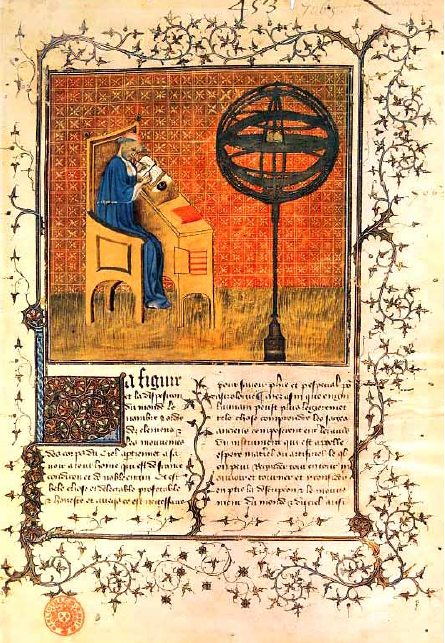}

\vspace{1cm}

\end{center}
\vfill

\newpage

\thispagestyle{empty}

\vspace*{15cm}

\noindent {\footnotesize \textsc{Stampato nel mese di aprile 2007\\
presso il Dipartimento di Matematica e Applicazioni,\\
Universit\`a degli Studi di Milano-Bicocca, via R. Cozzi 53, 20125
Milano, ITALIA}}.\\

\vspace*{.5cm}

\noindent {\footnotesize \textsc{Disponibile in formato elettronico
sul sito} \texttt{www.matapp.unimib.it}.\\
\textsc{Segreteria di redazione}: Ada Osmetti - Giuseppina Cogliandro\\
tel.: +39 02 6448 5755-5758  fax: +39 02 6448 5705}\\

\vspace*{.2cm} \noindent {\footnotesize \\ \textbf{ Esemplare fuori
commercio per il deposito legale agli effetti della Legge 15 aprile
2004 n.106}}.
\newpage
\begin{abstract}
We present here a detailed description of the model of ran-driven nuclear transduction in living cells
to be published in~\cite{Cangi09m}. The mathematical model presented is the first  to account for the active transport of molecules along the cytoplasmic microtubules. All parameters entering the models are thoroughly discussed. The simulations, carried out using the numerical scheme presented in~\cite{Cangi09DG},
reproduce the behavior observed experimentally.
\end{abstract}

\section{Introduction}


All cells receive and respond to signals from their surrounding.
Any external stimulus (ligand) acting on the cell {\em plasma membrane} activates several internal {\em second messenger} reactions that regulate virtually all aspects of cell behavior, including metabolism, movement, proliferation, and differentiation.
As living cells are highly compartmentalized systems in which biochemical reactions occur in physically distinct regions, the signal is transduced to the correct compartment to trigger the cellular response to the external environment.

The cell nucleus and specifically the genomic DNA contained in it, is the target of many intracellular transduction pathways. Indeed the response of the cell to the impinging signal
is obtained through the expression of specific genes.
In fact, as protein synthesis is carried out in the ribosome, the cellular response depends
as well on the export of RNAs out of the nucleus.

Here, we concentrate on molecular trafficking in and out the nuclear membrane. 
Molecules diffuse within the cytoplasm and reach the perinuclear region.
The translocation across the nucleocytoplasmic membrane or {\em nuclear envelope} (NE)
may proceed through the  {\em nuclear pore complexes} (NPCs) following two different mechanisms:
passive diffusion or facilitated translocation.
Only those molecules of mass less than 40 kDa~\cite{Gorlich03} can freely diffuse through the NPCs.
To permit the translocation of larger molecules, a system for active transport 
across the NPCs has evolved.
The cargo protein equipped with a {\em nuclear localization signal} (NLS) binds to a {\em nucleocytoplasmic transport receptor} (NTR) karyopherin known as 'importin', which mediates the transport through the NE.
%
The energy needed by facilitated translocation is provided by the Ran complex.
The Ran protein is a small GTPase (for a review, see for example: \cite{Quimby03})
cycling between two states: bound to guanosine triphosphate (RanGTP, active state) and to 
guanosine diphosphate (RanGDP, inactive state). The irreversible phosphorylation of RanGTP into RanGDP catalyzed by RanGAP in the cytoplasm maintains the Ran-cycle out of equilibrium permitting
the accumulation of RanGTP in the nucleus~\cite{Elbaum07}. This is then used to break the importin-cargo complex and thus permit the final release of cargo into the nucleus.

The mechanism of facilitated translocation permits the selective and regulated  translocation of relatively large molecules.
Indeed, it is understud that the combined action of importins, Ran complex and cargoes finely regulates the action of transcription factors within the nucleus.

Many mathematical models have been proposed that qualitatively reproduce the dynamics of intracellular signal transduction, often formulated in terms of ordinary differential equations describing the time course of the molecular concentrations (compartmental models)~\cite{Maca01,Smith02,Gorlich03,becs05,Riddick05,Khol06,Riddick07}. 
These are based on schematic descriptions of signal transduction pathways which do not take into account the spatial localization of the reactions~\cite{Harding05}. For a review on signaling networks modeling, see~\cite{Eungd04}.

In these pages we discuss spatial integrated models for Ran-driven nuclear import of molecules incorporating diffusion and membrane transport for a large-scale  model of living cell.
For the first time, we also include a model of the active transport along the {\em microtubules} of the importin-cargo complex.
The microtubules, together with the cytoplasmic filaments, constitute
the cytoplasm's dynamic structure that maintains cell shape (cytoskeleton).
The microtubules facilitates nuclear import of some proteins and viruses by providing a 
a preferential way directed towards the nuclear envelope (see~\cite{Roth07} and the references therein).
This is a typically spatial phenomena that cannot be taken into account in compartmental models.

To the best of our knowledge, Smith{\em et al.}~\cite{Smith02} is the only article presenting spatial simulations of cellular signal transduction pathways.
The scarcity of spatial mathematical models in the literature has to be related to the lack of knowledge of the relevant parameters, such as local diffusion and permeability coefficients. 
In~\cite{Smith02}, spatial and compartmental simulations are compared, giving similar answer. In our opinion, this has to be the case if we do not introduce any (spatial) details into the spatial model.
In fact, if the parameters are obtained by fitting experimental data with compartmental simulations, spatial models may as well be less realistic.

We shall discuss the crucial problem of parameters and pathways localization as we go through the various aspects of the model.

The Ran-driven nuclear import model 
discussed here is presented in~\cite{Cangi09m}.
The system of non-linear equations arising from spatial modeling shall be solved using
a new numerical technique based on Discontinuous Galerkin schemes. The details on the derivation
and numerical properties of the method are given in~\cite{Cangi09DG}. 

\section{The model}

The present model originates from the ODE model of Ran-driven nucleocytoplasmic import successively developed by Gorlich {\em et al.}~\cite{Gorlich01,Gorlich03}, Smith {\em et al.}~\cite{Smith02}, Riddick and Macara~\cite{Riddick05,Riddick07}, and Kopito and Elbaum~\cite{Elbaum07}.

Following~\cite{Elbaum07}, we keep the reaction network to the essentials.
We simplify the set of equations and explicitly introduce the spatial component within the variables which shall represent the molar concentrations of the molecular species.

Our model is a system of six semilinear parabolic PDEs set on two compartments: cytoplasm and nucleus.

In  the model we include, for each species, its initial concentration, molecular  reactions, Fickian diffusion, membranes translocation conditions, and facilitated transport through the microtubules.

In particular, transport through the NE and along the microtubules is only permitted to transport receptor complexes.

\subsection{Reaction network and mass action law}\label{sec:reactions}

The reaction terms are written in terms of the Law of Mass Action~\cite{deVries06}.
We assume that each biochemical reaction pathway  can be decomposed into unidirectional
elementary reaction and that  for each elementary reaction
the rate of change of the reactants is proportional to the product of their concentrations.
The Law of Mass Action is in fact a mathematical model expressing
the fact that the reaction depends on the number of molecular collisions
and the probability that a collision has enough kinetic energy to initiate the reaction.
The constant of proportionality is thus named {\em kinetic constant}.

Experimental values of the kinetic constants are usually available from {\em in vitro} experimentation from purified components, and thus they do not account for environmental interactions, competition and localizations. In fact, one of the goals of spatial modeling should be to obtain the correct `scaling relationship' of the parameters by fitting to known functional effects~\cite{Eungd04}.

Another limitation to quantitative mathematical modeling regards the rates of enzymatic irreversible reactions:
often the literature reports the Michaelis-Menten (MM) kinetic parameters  of the catalyzed reaction  instead of the kinetic constant of each reaction composing it
(see~\cite{MM13,Briggs25,SS89}, or the review in Chapter 6 of~\cite{MurrayI}).
Again, we have to stick to the data we are given, bearing in mind that an extra approximation has been introduced in the model.

\begin{figure}[h]
  \begin{center}
       \includegraphics[width=11.0cm]{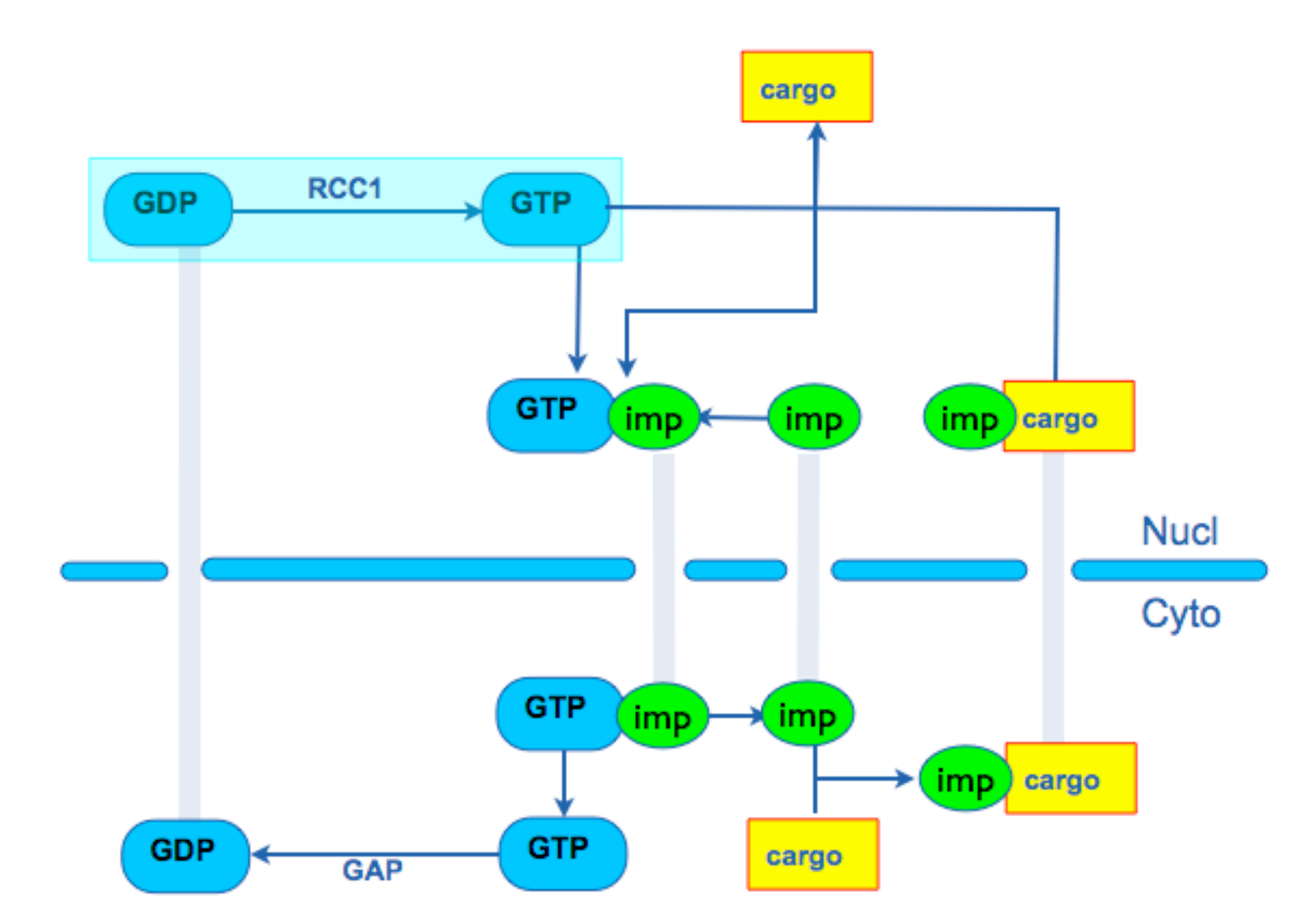}
  \end{center}
\caption {\footnotesize Reaction scheme of Importin-mediated transport between cytoplasm and nucleus.}
\label{fig:schema}
\end{figure}

The network of reactions involved in Ran-driven nuclear import pathway are schematize in Figure~\ref{fig:schema}  and described in Table~\ref{KineRan} using the species names of Table~\ref{tab:diff}.
\begin{table}[h]
{\small
  \begin{center}
     \begin{tabular}{|c|c|c|c|c|c|} \hline
\multicolumn{6}{|c|}{Reactions and kinetic constants for the  Ran transport model} \\ \hline
\tiny{Reaction} & \tiny{loc.} & \tiny{Term} & \tiny{const.} & \tiny{value (units)} & \tiny{ref.} \\ \hline
\tiny{$\Rt\stackrel{\Rg}{\rightarrow}\Rd$} & \tiny{cyto} & \tiny{$m_1(\Rt):=\displaystyle{\Kct\Rg\frac{\Rt}{\Kmt + \Rt}}$ }&
\tiny{$\begin{array}{c}
\vspace{1 mm}\Kct\\
\Kmt
\end{array}$}
 & 
\tiny{$\begin{array}{c}
\vspace{1 mm} 20.1\,\T{s}^{-1}\\
 0.7 \,\mu M 
\end{array}$}
 & \tiny{ \cite{Smith02}}
\\  \hline 
\tiny{$\Rd\stackrel{\Rc}{\rightarrow}\Rt$} & \tiny{nucl} & \tiny{$m_2(\Rd):=\displaystyle{\Kcd\Rc\frac{\Rd}{\Kmd + \Rd}}$ }&
\tiny{$\begin{array}{c}
\vspace{1 mm}\Kcd\\
\Kmd
\end{array}$}
 & 
\tiny{$\begin{array}{c}
\vspace{1 mm} 8.0 \,\T{s}^{-1}\\
 1.1 \,\mu M 
\end{array}$}
 &  \tiny{\cite{Gorlich03}}
\\  \hline 
\tiny{$\Rt + \Tr\rightleftharpoons\Tr_r $} & \tiny{both} & \tiny{
$\begin{array}{c}
\vspace{1 mm}
 r_1(\Rt,\Tr):= k_1 \, \Rt \Tr\\
r_{-1}(\Tr_r):= k_{-1}\, \Tr_r
\end{array}
$} &
\tiny{
$
\begin{array}{c}
\vspace{1 mm}
k_1\\
k_{-1}
\end{array}$} 
& 
\tiny{
$
\begin{array}{c}
\vspace{1 mm}
0.1\, (\mu M\:\T{s})^{-1}\\
0.3\, \T{s}^{-1}
\end{array} $}   & \tiny{\cite{Riddick05}}
\\ \hline 
\tiny{$\cargo +\Tr\rightarrow \Tr_c$} & \tiny{cyto} & \tiny{$ r_2(\cargo,\Tr):= k_2 \, \cargo \,\Tr$} &
\tiny{$k_{2}$} 
& 
\tiny{$0.15\, (\mu M\:\T{s})^{-1}$}   & \tiny{\cite{Riddick05}}
\\ \hline 
\tiny{$\Rt+\Tr_c\rightarrow \Tr_r+\cargo$} & \tiny{nucl} & \tiny{$ r_3(\Rt,\Tr_c):= k_3 \, \Rt\,\Tr_c$} &
\tiny{$k_{3}$ }
&
\tiny{$0.1\, (\mu M\:\T{s})^{-1}$}   & \tiny{\cite{Riddick05}}
\\ \hline 
  \end{tabular}
 \end{center}}
\caption{\footnotesize Kinetic constants for the reactions involved in the Ran-driven transport process,
with respect to the species names in Table~\ref{tab:diff}. Second order constants are measured in $\mu$M$^{-1}$ $\T{s}^{-1}$ while first order constants are measured in $\T{s}^{-1}$. RCC1 catalyzed reaction is represented with a MM scheme in terms
of the constant concentration $c1$, as in Smith e al, 2002. For a more refined scheme 
(a multistep scheme) see G\"orlich et al, 2003. The MM scheme gives a good approximation as 
 the intermediate steps (complex formation) are more rapid than the exchange reaction. 
The enzymes concentrations are assumed to remain constant. In particular, following\cite{Riddick05},
we set $\Rg=0.5 \mu$M and $\Rc=0.7 \mu$M.
}
\label{KineRan}
\end{table}

We shall assume that facilitated translocation through the nucleocytoplasmic membrane is only permitted to molecules associated to  a transport receptor (see below the description on membrane shuttling).

Receptors (adapters)-mediated import necessitates  RanGTP to disassociate the cargo from the receptor ones the complex has entered the nucleus.
Thus, the transport mechanism relies on the presence of large concentrations of RanGTP in the nucleus in comparison  to the cytoplasm (Ran gradient).
The Ran complex is responsible to maintain the RanGTP/RanGDP gradient, thus permitting cargo accumulation in the nucleus.

Within the cell, the small GTPase Ran can bind to guanosine nucleotide phosphates GDP and GTP, forming the RanGDP and RanGTP complexes, respectively. These reaction, described in the model by MM-kinetics, are catalyzed by two specific enzymes: RanGAP which is located in the cytoplasm and RanGEF (RCC1) located within the nucleus. 
Due to this cycle of reactions,
a large concentration jump of RanGTP across the NE is readily established, with a high RanGTP concentration in the nucleus.

The first cycle of cargo import begins with the formation of the cargo receptor complex.
The receptor involved in nuclear import is made of two subunits:
the {\em adaptor} Importin-$\alpha$  which binds the cargo Nuclear Localization Signal (NLS) and the {\em receptor} Importin-$\beta$ which permits the translocation through the NPC. Here, we simplify matters by considering the binding of the cargo with the whole importin complex (see~\cite{Riddick05,Riddick07} for a detailed analysis of the various import pathways).

The second (parallel) cycle is the nuclear import of
cytoplasmic RanGDP by the effector NTF2. Here we assume that NTF2 is always available and treat  RanGDP as all bound to NTF2.
Nuclear RanGDP, imported by the transport receptor NTF2, interacts with the nuclear RanGEF (RCC1) which catalyzes the exchange between nucleotide GDP and GTP forming the nuclear RanGTP. The complete biochemical scheme of this exchange is composed of a sequential set of enzymatic steps~\cite{Klebe95a,Klebe95b,Gorlich03,Riddick05} which, following~\cite{Smith02},  we schematize as a single reaction.

The nuclear importin-cargo complex can bind the RanGTP to form two complexes, RanGTP$\cdot$Imp$\beta$ and Imp$\alpha$$\cdot$cargo. This latter complex can dissociate, create free Imp$\alpha$ and cargo or interact with RanGTP to form the RanGTP$\cdot$Imp$\alpha$ complex and free cargo.
We simplify this pathway with the single reaction
\[
T\bullet C + R \rightarrow C + T \bullet R,
\]
written in terms of  the generic transport receptor $T$, cargo $C$, and with $R$ representing the concentration of Ran (bound to GTP in this case).

The free cargo, which usually is an activated transcription factor, can now bind DNA and activate the gene expression program. In our model, these reactions are neglected, and thus the cargo is allowed to accumulate within the nucleus,
while the receptor-RanGTP complex exits the nucleus and, eventually, dissociates.


\subsection{Diffusion}

The molecular species diffusion is expressed in terms of Fick's law. All the molecular species diffuse with a specific diffusion coefficient.
This is obtained from the formula
\[
D=\frac{KT}{6\pi \eta R_s},
\]
in terms of the Boltzmann constant $K$, the temperature $T$, 
the Stokes radius $R_s$, and the viscosity of the medium $\eta$ (see {\em e.g.}~\cite{Gorlich01}).

The viscosity of the cytoplasm has been measured in the nineties for a wide range of molecular weights and cellular environments~\cite{Seksek97}.
It was initially thought that the viscosity had to depend on the molecular size to account for
the sieving effect of the cytoskeletal network and other macromolecular structures.
For this reason, solutes diffusion was described by
the {\em translational diffusion coefficient}~\cite{Kao93}, depending on
the viscosity of the  fluid-phase cytoplasm and on collisions with intracellular components.
The former, representing the viscosity sensed by a small solute in the absence of interactions with macromolecules and organelles,  is $1.2$--$1.4$ times that of water~\cite{Verkman91}.

The net viscosity of the cytoplasm of Swiss 3T3 fibroblasts with respect to small solutes such as metabolites and nucleic acids was measured by Kao {\em et al.}~\cite{Kao93}  using fluorescence recovery after spot-photobleaching (FRAP)~\cite{Axelrod76}. They found that mobility was around four to fivefold slowed than in water.
Contrary to expectations, later studies extended the validity of this result to the 
mobility of macro-molecules of up to few hundreds of kDs of mass~\cite{Seksek97}.
Moreover, the same applies to the nucleus~\cite{Seksek97} and the mitochondrial matrix\cite{Parti98}, thus viscosity fourfold higher than water can be assumed throughout.

A correction may be employed in the proximity of the plasma membrane as the 
translational diffusion of the fluorescent probe BCECF already used by~\cite{Kao93} 
is found to be  twofold lower near the cell membrane  due to high density of proteins~\cite{Swamin96}.



The diffusion coefficients used in the model are reported in Table~\ref{tab:diff}.
\begin{table}[h]
  \begin{center}
     \begin{tabular}{|c|c|c|c|c|} \hline
\multicolumn{5}{|c|}{Diffusion coefficients for Ran model} \\ \hline
reactant & variable     &  diff. coef.    &   value ($\mu$m$^2$/s) & reference \\ \hline
RanGTP   & $\Rt$        & $d_r$           & 22             & \cite{Caudron05}  \\ \hline
RanGDP (with NTF2)& $\Rd$ & $d_r$         & 20             & \cite{Gorlich03,Smith02}  \\ \hline
cargo    & $\cargo$     & $d_c$           & 18.2           & calculated \\ \hline
receptor & $\Tr$        & $d_t$           & 14             & \cite{Caudron05} \\ \hline
RanGTP$\bullet$receptor & $\Tr_r$ & $d_{tr}$   & 14        & \cite{Caudron05}  \\ \hline
 cargo$\bullet$receptor & $\Tr_c$ & $d_{tc}$   & 12.4      & calculated \\ \hline 
  \end{tabular}
 \end{center}
\caption{\footnotesize 
Diffusion coefficients for the molecular species involved in the Ran-driven
 transport system. Values are the same for the cytoplasm and the nucleus. 
These data have been calculated assuming a cytoplasm viscosity 5 times higher than water~\cite{Seksek97}   
at a temperature of $20^{\circ}$ Celsius.
Cargo weight taken as the ERK mass, 42 kDa (for ERK1) or 44 kDa (ERK2).
Data from the literature are not always consistent: for instance, Smith {\em et al}~\cite{Smith02}
estimates the diffusivity of Ran in $30 \mu$m$^2$/s.
}
\label{tab:diff}
\end{table}

\subsection{Facilitated translocation through the nuclear envelope}

In this report, we concentrate on nuclear import by assuming that the  cargo is already in the cell  periphery and that there is no exchange of substances between the cell and  the surrounding environment through the {\em plasma membrane}.
Mathematically, no-flux conditions are specified for all species at the plasma  membrane. Thus, given the generic species concentration $u$ we impose:
\begin{equation}\label{eq:noflux}
d^c_u\frac{\partial{u}^c}{\partial\mbf{n}}=0,
\end{equation}
where $\mbf{n}$ is the unit normal vector pointing outside the cell.
In fact, the mechanism of signal transduction through the plasma membrane is often quite different from that of NE trapassing. The external signals are transduced from external membrane receptors
into internal second messenger activation on the inside of the plasma membrane~\cite{Alberts}. Thus, no passage of matter through the membrane is involved.

All shuttling across the nuclear envelope takes place through the nuclear pore complexes (NPCs), multi-component protein structures spanning the nuclear envelope.  Two mechanisms of translocation are permitted:
passive diffusion or facilitated translocation.
Due to the limited diameter of the pore lumen (about 10 nm),
ions and small metabolites can freely diffuse through the NPCs only if their mass is less than $\sim \hspace{-1 mm}40$ kDa.
Molecules of size greater than that of the NPCs can still shuttle across the NE
by facilitated translocation.
The actual mechanism of facilitated translocation is still under research (for a review of various models, see~\cite{Aebi03}).
The molecule must posses a specific aminoacid domain which binds to specific transport receptors (importins and exportins) that operate a conformational change in the NPC whose `functional size' is in the order of 25 nm~\cite{Mattaj98,becs05}.

Thus, we assume that facilitated translocation is allowed only to receptor complexes~\cite{Elbaum07}.

The study presented in~\cite{Elbaum07} proves that the accumulation of cargo in the nucleus
is not due, as previously thought, to the ability of the `importin' receptor to cross the NE is
unidirectional (from the cytoplasm to the nucleus).
The translocation is actually bidirectional, with the accumulation of cargo in the nucleus
being rather due to the asymmetric concentration of RanGTP. 

In the absence of a detailed and realistic model of the reactions happening within the pore
(interactions of proteins with nucleoporins, meshwork of filaments within the lumen of the pore, etc), a more convenient approach is to use a `coarse-grain' formulation, in terms of permeabilities times 
concentration differences between the two compartments (see~\cite{Smith02}).
The flux across the nucleo-cytoplasmic boundary is modeled as the product of a proportionality
constant (the permeability $P$) times the concentration difference across the nucleocytoplasmic 
boundary, as a direct consequence of the fact that the  flux does not require additional energy input \cite{Smith02}. Thus, 
we fix the following
Kedem-Katchalsky conditions at the nuclear envelope:
\begin{equation}
\label{permeq}
d_u^n\frac{\partial u^n}{\partial \mbf{n}}=p_u({u}^c-{u}^n),
\end{equation}
where $\mbf{n}$ is the unit vector pointing from the nucleus to the cytoplasm and the factor ${u}^n-{u}^c$ is the concentration difference of species $u$ across the nucleocytoplasmic boundary.
The  nuclear membrane permeability is assumed constant as experimental data indicate binding to the pore complex is far from saturation \cite{Gorlich01,Smith02}.
The permeability values can be calculated from experimental values on the capacity,area, and number  of the NPCs present on the nuclear envelope.
Ribbeck and Gorlich~\cite{Gorlich01}  estimated $>100$ translocation events/NPC/second and 
a density of NPCs of $5.1 \pm 0.2 $ NPCs$/\mu m^2$ in the NE of HeLa cells.

Estimates of  translocation rates of most molecular species are available from the literature and have already been used in a number of studies~\cite{Maca01,Smith02,Gorlich03,Riddick05,Riddick07}.
The values used in the simulations are taken from~\cite{Smith02} and shown in Table~\ref{tab:perm}.
These were fitted by comparison of experimental data with compartmental simulations: more investigation on the correct values for spatial
modeling is needed, for instance by following the approach of
the homogenized `effective' permeabilities proposed in~\cite{Berez04,Berez06}).
\begin{table}[h]
{\small
  \begin{center}
     \begin{tabular}{|c|c|c|c|} \hline
\multicolumn{4}{|c|}{NPC permeabilities for Ran model}       \\ \hline
reactant & perm. const. & value ($\mu$m/s) & reference(s)             \\ \hline
$\Rd$ (with NTF2)   &$p_d$  & 3.73      & \cite{Smith02}   \\ \hline
$\Tr$               &$p_t$  & 1.87     & \cite{Smith02} \\ \hline
$\Tr_r$             &$p_{tr}$ & 1.87      & \cite{Smith02} \\ \hline
$\Tr_c$             &$p_{tc}$  & 1.87      & \cite{Smith02} \\ \hline 
  \end{tabular}
 \end{center}}
\caption{\footnotesize Permeability values for molecular species involved in membrane translocation. The permeability is the kinetic effect of the molecular species with the nucleoporin complex (NPC).}\label{tab:perm}
\end{table}

The set of boundary conditions is closed by imposing continuity of the flux:
\begin{equation}
\label{contflux}
d^c_u\frac{\partial{u}^c}{\partial\mbf{n}}=d^n_u\frac{\partial{u}^n}{\partial\mbf{n}},
\end{equation}

For molecular species which cannot pass through the nuclear membrane we just impose homogeneous
Neumann boundary conditions on the relevant compartment. For instance, if the species $u$ is
confined in the cytoplasm, then
\begin{equation}
d^c_u\frac{\partial{u}^c}{\partial\mbf{n}}=0.
\end{equation}


\subsection{Active transport along the microtubules}

Microtubules are known to be used by macro molecules like
adenoviruses~\cite{Dinh06} and by intracellular organelles like endocytic vesicles~\cite{Dinh06} as preferential ways of motion. 
The adenoviruses are $~90 nm$ in diameter 
and thus their diffusion speed in the cytoplasm is very limited: 
for this reason they must resort in  active transport along the microtubules in order to reach the nucleus. 
%
%

It is by now well established that also some smaller molecules
utilize motor-assisted transport along microtubules (MTs)~\cite{Gunder99,Camp03,ElbaumMT05}.
Examples goes from mRNA~\cite{Tekotte02} localizing after
nuclear export, to the tumor suppressor protein
p53~\cite{Gianna02} and Rb~\cite{Roth07}.

Many similarly sized molecules resort on passive diffusion only to reach the nucleus NE.
Thus, in this respect, active transport is not essential to the transduction precess.
Rather it must be seen as a way  to improve the transduction efficiency.
Our aim is to introduce a complete model of the transduction mechanism, accounting for both diffusion and active transport along the micrutubules, that can be used to better understand the characteristics and importance of the latter mechanism.

The microtubules,  which, together with the cytoplasmic filaments, constitute the cytoskeleton, are usually organized into a single array centered near the nucleus. They are characterized by an orientation, with the minus ends located near the NE.
Active transport along the MTs is permitted by binding to a {\em motor} protein~\cite{Nedel01,Smith01}.
We distinguish among two families of motor proteins: dynein, which permits transport from the plus end to the
minus end, and kinesin, which transports in the opposite direction.
For instance, kinesin anterograde transport  along microtubules is known to be used by particles for covering large distances along nerve axons (see~\cite{Smith01} and the references therein).

The receptor Importin-$\alpha$ is able to bind to the MTs via motor proteins, thus permitting the active transport of NLS cargoes towards the nucleus.
It is this phenomenon that we wish to include in our model.
The speed of transport of a motor protein attached to an MT is of about $0.5$ to $1\, \mu m\, \T{s}^{-1}$~\cite{Smith01,Nedel01},
but pauses and changes of direction of motion are often observed,
suggesting that molecules may proceed by detaching and changing type of motor~\cite{Smith01}.
Tentative mathematical models to fully describe MTs effective directional 
transport have been proposed, for instance, 
in the References~\cite{Smith01,Nedel01,Dinh05,ElbaumMT05,ElbaumMT08}.

Here, we assume that the given cargo-receptor complex can only bind to dynein, and thus MT transport is directed from the cell periphery
to the nuclear envelope, and verify if the model reproduces the experimental
results of Roth {\em et al.}~\cite{Roth07}.
We simulate active transport of Nuclear Localization Signals (NLS) along MTs
by introducing an advection term in the species flux, which becomes:
\[
J_u=d_u \nabla u -{\mathbf{b}} u.
\]
The advection field $\mathbf{b}$ is taken  of constant modulus equal to $1$ pointing towards the nucleus
(unidirectional transport from the plus end to the minus end motored by dynein).

\subsection{Partial differential equations model}

Let $\Omega$ represent the cell's domain. Further, we denote by   $\partial\Omega$ its boundary (the plasma membrane), and by $\Gatr$
the interface between cytoplasm and nucleus (the nuclear envelope).

On both membranes, we fix the normal unit vector
$\normale$ as  pointing inside the cytoplasm.
Accordingly, the jump on $\Gatr$ of a given species concentration $u$ will be denoted by
\[
\jumptwo{u}:=u^c|_{\Gatr}-u^n|_{\Gatr}.
\]

By collecting all the model's term described above, we obtain the following
system of semilinear parabolic equations.
In the cytoplasm, the species concentrations obey:
\[
\left\{
\begin{array}{lll}
\vmm{2}
\displaystyle\frac{\partial{\Rt}}{\partial{t}} \hmm{-2}&=& \hmm{-2} 
d_{r}\Delta\Rt
 - m_1(\Rt)-r_1(\Rt,\Tr)+ r_{-1}(\Tr_r), 
\\
\vmm{2}
\displaystyle\frac{\partial{\Rd}}{\partial{t}} \hmm{-2}&=& \hmm{-2}
d_{r}\Delta\Rd
+ m_1(\Rt),
\\
\vmm{2}
\displaystyle\frac{\partial{\Tr_r}}{\partial{t}} \hmm{-2}&=&\hmm{-2} d_{tr}\Delta\Tr_r
+ r_1(\Rt,\Tr)-r_{-1}(\Tr_r),
\\
\vmm{2}
\displaystyle\frac{\partial{\cargo}}{\partial{t}}\hmm{-2} &=&\hmm{-2} d_{c}\Delta\cargo
-r_2(\cargo,\Tr),
\\
\vmm{2}
\displaystyle\frac{\partial{\Tr}}{\partial{t}} \hmm{-2}&=&\hmm{-2}
d_{t}\Delta\Tr
 -r_1(\Rt,\Tr)+ r_{-1}(\Tr_r)-r_2(\cargo,\Tr),
\\
\displaystyle\frac{\partial{\Tr_c}}{\partial{t}}\hmm{-2}&=&\hmm{-2}
d_{tc}\Delta\Tr_c
-\nabla({\mathbf{b}}\,\Tr_c)
+ r_2(\cargo,\Tr),
\end{array}
\right.
\]
Notice the advective term modeling transport of the receptor-cargo complex which
appears in the last equation above.

As we assume that no matter is entering or exiting the cell through the plasma
membrane, the homogeneous Neumann boundary condition~(\ref{eq:noflux}) is
imposed to all species.

In the nucleus, we have the following system of reaction-diffusion equations:
\[
\left\{
\begin{array}{lll}
\vmm{2}
\displaystyle\frac{\partial{\Rt}}{\partial{t}} \hmm{-2}&=& \hmm{-2} 
d_{r}\Delta\Rt
+ m_2(\Rd)-r_1(\Rt,\Tr)+ r_{-1}(\Tr_r)
-r_3(\Rt,\Tr_c),
\\
\vmm{2}
\displaystyle\frac{\partial{\Rd}}{\partial{t}} \hmm{-2}&=& \hmm{-2}
d_{r}\Delta\Rd
- m_2(\Rd),
\\
\vmm{2}
\displaystyle\frac{\partial{\Tr_r}}{\partial{t}} \hmm{-2}&=&\hmm{-2} d_{tr}\Delta\Tr_r
+ r_1(\Rt,\Tr)-r_{-1}(\Tr_r)+r_3(\Rt,\Tr_c),
\\
\vmm{2}
\displaystyle\frac{\partial{\cargo}}{\partial{t}}\hmm{-2} &=&\hmm{-2} d_{c}\Delta\cargo
+r_3(\Rt,\Tr_c),
\\
\vmm{2}
\displaystyle\frac{\partial{\Tr}}{\partial{t}} \hmm{-2}&=&\hmm{-2}
d_{t}\Delta\Tr
 -r_1(\Rt,\Tr)+ r_{-1}(\Tr_r),
\\
\displaystyle\frac{\partial{\Tr_c}}{\partial{t}}\hmm{-2}&=&\hmm{-2}
d_{tc}\Delta\Tr_c
+\nabla({\mathbf{b}}\,\Tr_c)
-r_3(\Rt,\Tr_c),
\end{array}
\right.
\]

The two systems of equations above are coupled through the appropriate transmission conditions on $\Gatr$:
\[
\left\{
\begin{array}{lll}
\vmm{2}
d_r\frac{\partial{\Rt^{c,n}}}{\partial\normale}|_{\Gatr}\hmm{-3.5}&=0,&
\\
\vmm{2}
d_r\frac{\partial{\Rd^c}}{\partial\normale}|_{\Gatr}\hmm{-3.5}&=
d_r\frac{\partial{\Rd^n}}{\partial\normale}|_{\Gatr}\hmm{-4}&=
p_d\jumptwo{\Rd},
\\
\vmm{2}
d_{tr}\frac{\partial{\Tr_r^c}}{\partial\normale}|_{\Gatr}\hmm{-3.5}&=
d_{tr}\frac{\partial{\Tr_r^n}}{\partial\normale}|_{\Gatr}\hmm{-4}&=
p_{tr}\jumptwo{\Tr_r},
\\
\vmm{2}
d_c\frac{\partial{\cargo}^{c,n}}{\partial\normale}|_{\Gatr}\hmm{-3.5}&=0,&
\\
\vmm{2}
d_{t}\frac{\partial{\Tr^c}}{\partial\normale}|_{\Gatr}\hmm{-3.5}&=
d_t\frac{\partial{\Tr^n}}{\partial\normale}|_{\Gatr}\hmm{-4}&=
p_{t}\jumptwo{\Tr},
\\
d_{tc}\frac{\partial{\Tr_c^c}}{\partial\normale}|_{\Gatr}\hmm{-3.5}&=
d_{tc}\frac{\partial{\Tr_c^n}}{\partial\normale}|_{\Gatr}\hmm{-4}&=
p_{tc}\jumptwo{\Tr_c}.
\end{array}
\right.
\]

These equations needs to be provided with initial conditions which are specified below in Section~\ref{sec:numexp}.

\section{Numerical approximation of the mathematical model}

We name the two subdomains of $\Omega$ representing the cytoplasm and the nucleus as
$\Omega_1$  and $\Omega_2$, respectively. Thus, $\Omega=\Omega_1\cup\Omega_2\cup\Gatr$, where $\Gatr:=\bar{\Omega}_1\cap\bar{\Omega}_2$.

In order to ease notation in the presentation of the numerical method, we rewrite  the mathematical
model described in the previous section in vector form.

\subsection{The pde problem in vector form}

Let us denote the number of unknown concentrations by $n$, 
and collect all concentrations variables in  the vector function
\[
\bu:=(u_1,\dots,u_n)^T:\Omega_1\cup\Omega_2\to\mathbb{R}^n.
\]

We shall look for solutions belonging to the space
$\honen:=[H^1(\Omega_1\cup\Omega_2)]^{n}$ at any time $t\in[0,T]$, with $T$ representing
some final time.

We introduce the diagonal matrix  $U=\diag(u_1,\dots,u_n)$ and the gradient $\nabla \bu:\Omega_1\cup\Omega_2\to\mathbb{R}^{n\times d}$
given by $\nabla \bu:=(\nabla u_1,\dots,\nabla u_n)^T$,
with $\nabla u_i(\mbf{x})\in\mathbb{R}^{d}$, $i=1,\dots, n$, $\mbf{x}\in\Omega_1\cup\Omega_2$.
Further, for a tensor $\mbf{Q}:\Omega_1\cup\Omega_2\to\mathbb{R}^{n\times d}$, with rows $Q_i$, $i=1,\dots, n$,
we define its divergence $\nabla \cdot\mbf{Q}:\Omega_1\cup\Omega_2\to\mathbb{R}^{n}$
by $\nabla \cdot\mbf{Q}:=(\nabla \cdot Q_1,\dots,\nabla \cdot Q_n)^T$.

We collect all reaction terms in the vector field 
$\mbf{f}(\bu)$.
Similarly, the advection and diffusion coefficients are
collected in the tensors 
$\mbf{B}\in [C^1(0,T;\bar{\Omega}\backslash\Gatr)]^{n\times d}$, whose rows are denoted by $B_i$, $i=1,\dots, n$, and
$A\in [C(0,T;\Omega_1\cup\Omega_2)]^{n\times n}$ diagonal, with $A=\diag(d_1,d_2,\dots,d_n)$,
where $d_i>0$, $i=1,\dots, n$
are the diffusivity of the various species.
Notice that $B_i$ will be the zero vector if the $i$-th
species is not transported. In particular, we assume that
every $B_i$ is continuous and supported on an open (eventually empty) subset of $\Omega_1$ (the MTs are located in the cytoplasm).
The implication of this assumption which is used below in the definition of the method is that
$\mbf{B}=0$ on the boundaries of both $\Omega_1$ and $\Omega_2$.
Finally, let $\bu_0\in [L^2(\Omega)]^n$ represent the initial conditions.

We look for solutions $\bu\in L^2(0,T;\honen)$ of the following initial and boundary value problem:
\begin{equation}\label{eq:theproblem}
\left\{
\begin{array}{lll}
\vmm{2}
\bu_t = \nabla\cdot (A\nabla \bu - 
U\mbf{B})+\mbf{f}(\bu) &\text{in }\: [0,T]\times(\Omega_1\cup\Omega_2),
\\
\vmm{2}
\bu(0,x)=\bu_0(x)&\text{on }\:\{0\}\times\Omega,
\\
\vmm{2}
\big(A\nabla \bu\big)\mbf{n} = \mbf{0} &\text{on }\: \partial\Omega,
\\
\vmm{2} 
\big(A\nabla \bu\big)\mbf{n}|_{\Omega_1} =
{\rm P}\: (\bu|_{\Omega_2}-\bu|_{\Omega_1})&\text{on }\: \Gatr,\\
\big(A\nabla \bu\big)\mbf{n}|_{\Omega_2} =
{\rm P}\: (\bu|_{\Omega_1}-\bu|_{\Omega_2})&\text{on }\: \Gatr,
\end{array}
\right.
\end{equation}
where $\mbf{n}|_{\Omega_j}$ denotes the unit outward normal vector to $\Gatr$ from $\Omega_j$, $j=1,2$,
and 
${\rm P}:=\diag(p_1,\dots,p_n)$ the diagonal matrix of permeabilities $p_i$.

More general boundary conditions, including mixed-type conditions,  and transmission conditions, including nonlinear conditions, are discussed in~\cite{Cangi09DG}.

\subsection{The interior penalty discontinuous Galerkin method}

In recent years, there has been an increasing interest in \emph{discontinuous
Galerkin finite element methods} (DGFEM)\cite{cks}. The renewed interest in these
methods is due to their very good stability properties when used to approximate solutions to convection-dominated convection-diffusion problems~\cite{CockburnNASAVKI99}, as well as due to the great flexibility in grid design they offer. Moreover, they naturally embed
good local conservation properties of the state variable, which can
be advantageous for time-dependent problems.

Regarding the application of DGFEMs to the present problem, we add that the weak imposition of
boundary conditions typical of discontinuous methods
permits a very natural imposition of the transmission conditions.
The good stability properties of the method on convection-dominated problems are relevant to
the case of cargos representing large viruses, not considered
in the simulations below.
Indeed, the diffusion of macromolecular structures is impeded
by the cytoskeleton (the effective viscosity is up to
a hundred times higher than that of water),   and
active transport along the MTs becomes the dominating mean of migration towards the nucleus~\cite{Camp03}.

The method used here,  which is detailed in~\cite{Cangi09DG}, is a generalization of the symmetric version of the interior penalty discontinuous Galerkin (IPDG) method. This was introduced in~\cite{houston-etal:01} for the 
solution of second-order partial differential equations with nonnegative characteristic form
and later extended to the solution of semilinear time-dependent problems in~\cite{Suli07}.


The finite element method is introduced by defining a shape-regular subdivision
$\mathcal{T}$ of $\Omega_1\cup\Omega_2$ into disjoint open {\em elements} $\el$.
We assume that the  elements of $\mathcal{T}$ `belong' to either subdomain, so that
subdivisions of  $\Omega_1$ and $\Omega_2$  are authomatically introduced.
Abusing notation, we use the simbol $\Gatr$ to denote the union of elemental faces
belonging to both subdivisions.
Further we decompose the subdivision {\em skeleton} $\Gamma:=\cup_{\el}\partial\kappa$ into three disjoint subsets
\[
\Gamma=\partial\Omega\cup\Gamma_{\dint}\cup\Gatr,
\]
where $\Gamma_{\dint}:=\Gamma\backslash(\partial\Omega\cup\Gatr)$.

We assume that the subdivision $\mathcal{T}$
is constructed via mappings $F_{\kappa}$, where
$F_{\kappa}:\hat{\kappa}\to\kappa$ are smooth maps
with non-singular Jacobian, and $\hat{\kappa}$ is the reference $d$-dimensional simplex or the reference $d$-dimensional (hyper)cube; the
maps are assumed to be constructed so as to ensure that the
union of the closures of the elements $\el$ forms a covering of the
closure of $\Omega$, i.e., $\bar{\Omega}=\cup_{\el}\bar{\kappa}$.

For a nonnegative integer $m$, we denote by $\mathcal{P}_m(\hat{\kappa})$, the set of polynomials of total degree at most $m$ if $\hat{\kappa}$ is the
reference simplex, or the set of all tensor-product polynomials on $\hat{\kappa}$ of
degree $k$ in each variable, if $\hat{\kappa}$ is the reference hypercube. We consider
the $hp$-discontinuous finite element space
\begin{equation}\label{eq:FEM-spc}
S:=\{v\in L^2(\Omega):v|_{\kappa}\circ F_{\kappa}
    \in \mathcal{P}_{m_{\kappa}}(\hat{\kappa}),\,\el\}.
\end{equation}

Next, we introduce some trace operators. Let $\kappa^+$, $\kappa^-$ be two (generic) elements sharing an edge
$e:=\partial\kappa^+\cap\partial\kappa^-\subset\Gamma_{\dint}\cup \Gatr$.
Define the outward normal unit vectors $\mbf{n}^+$ and $\mbf{n}^-$ on $e$ corresponding
to $\partial\kappa^+$ and $\partial\kappa^-$, respectively. For functions
$\mbf{q}:\Omega\to\mathbb{R}^n$ and $\mbf{Q}:\Omega\to\mathbb{R}^{n\times d}$ that may be discontinuous across $\Gamma$,
we define the following quantities.
For
$\mbf{q}^+:=\mbf{q}|_{e\subset\partial\kappa^+}$, $\mbf{q}^-:=\mbf{q}|_{e\subset\partial\kappa^-}$, and
$\mbf{Q}^+:=\mbf{Q}|_{e\subset\partial\kappa^+}$, $\mbf{Q}^-:=\mbf{Q}|_{e\subset\partial\kappa^-}$, we set
\[ \mean{\mbf{q}}:=\ha(\mbf{q}^+ + \mbf{q}^-),\ \mean{\mbf{Q}}:=\ha(\mbf{Q}^+ + \mbf{Q}^-),\]
and
\[ \jumpthree{\mbf{q}}:=\mbf{q}^+\otimes \mbf{n}^++\mbf{q}^-\otimes \mbf{n}^-,
\ \jumptwo{\mbf{Q}}:=\mbf{Q}^+ \mbf{n}^++\mbf{Q}^- \mbf{n}^-,\]
where $\otimes$ denotes the standard tensor product operator, whereby $\mbf{q}\otimes\mbf{w}=\mbf{q}\mbf{w}^T$.
If $e\in \partial\kappa\cap\Gamma_{\partial}$, these definitions are modified as follows
\[ \mean{\mbf{q}}:=\mbf{q}^+ ,\ \mean{\mbf{Q}}:=\mbf{Q}^+,
\  \jumpthree{\mbf{q}}:=q^+\otimes\mbf{n},\ \jumptwo{\mbf{Q}}:=\mbf{Q}^+\mbf{n}.\]

Finally, we introduce the mesh quantities
$\tth:\Omega\to \mathbb{R}$, ${\tt m}:\Omega\to \mathbb{R}$, by $\tth(x)= \diam{\kappa}$,
 ${\tt m}(x)= m_{\kappa}$, if $x\in \kappa$, and $\tth(x)= \mean{\tth}$, if $x\in \Gamma$, ${\tt m}(x)= \mean{{\tt m}}$, if $x\in \Gamma$, respectively.

In order to define the IPDG finite element method, we further introduce 
$\Sigma:=\diag(\sigma_1,\dots,\sigma_n)$ the diagonal matrix containing the
\emph{discontinuity-penalisation parameters} $\sigma_i:\Gamma \backslash \Gatr$,
and denote $\mathcal{B}:=1/2\diag(|B_1\cdot \mbf{n}|,\dots,|B_n\cdot \mbf{n}|)$.

The IPDG-in-space method for the system (\ref{eq:theproblem})
reads: 
find $\bu_h\in L^2(0,T;[S]^n)$ such that
\begin{equation}\label{eq:semidiscrete}
\ltwoin{(\bu_h)_t}{\bv_h}+B(\bu_h,\bv_h)= \ltwoin{\mbf{f}(\bu_h)}{\bv_h},
\quad\texte{for all} \bv_h\in L^2(0,T;[S]^n),
\end{equation}
where $\ltwoin{\cdot}{\cdot}$ is the $L^2$ scalar product, the bilinear form
 and $B(\cdot,\cdot)$ is defined by
\begin{equation}\label{eq:bilinear}
\begin{aligned}
B(\bu_h,\bv_h):=&\su\int_{\kappa} (A\nabla\bu_h-U_h\mbf{B}): \nabla \bv_h
+\int_{\Gatr} {\rm P}\jumpthree{\bu_h}:\jumpthree{\bv_h}\\
&\hmm{- 15}-\int_{\Gint}\Big(\mean{A\nabla \bu_h-U_h\mbf{B}}:\jumpthree{\bv_h}+\mean{A\nabla \bv_h}:\jumpthree{\bu_h}
-(\Sigma+\mathcal{B})\jumpthree{\bu_h}:\jumpthree{\bv_h}\Big).
\end{aligned}
\end{equation}
A fully discrete formulation is obtained from~(\ref{eq:semidiscrete}) by using any standard 
time stepping method like, for instance any Runge-Kutta time stepping. Full details on the derivation
and analysis of the method can be found in~\cite{Cangi09DG}.

\section{Numerical simulations}\label{sec:numexp}

We use the numerical method described above to simulate our model of cargo import. 
\begin{table}[h]
  \begin{center}
     \begin{tabular}{|c|c|c|c|} \hline
\multicolumn{4}{|c|}{Initial concentrations} \\ \hline
reactant & localization & $\mu$M  & reference(s) \\ \hline
$\Rt$ & cyto & 3     & \cite{Riddick05,Smith02}          \\ \hline
$\Rd$ & cyto & 3     & \cite{Riddick05,Smith02}          \\ \hline
$\cargo$ & near pl. membrane & 8        & \cite{Riddick05}          \\ \hline
$\Tr$ & cyto & 4       & \cite{Riddick05}            \\ \hline
  \end{tabular}
 \end{center}
\caption{\footnotesize Initial concentrations of constituents 
of the Ran system used in the model by Riddick {\em et al.}~\cite{Riddick05}. All other initial concentrations are set to zero. The concentrations of the enzymes $RCC1$ and RanGAP are assumed to be constant, as described in Section~\ref{sec:reactions}}\label{tab:conc}
\end{table}
The goal is to verify the experimental results 
presented by Roth {\em et al.}~\cite{Roth07}. 
Two experiments are presented in~\cite{Roth07}: in the first, the NLS cargo known to
bind to the MTs is activated in a cell with intact cytoskeleton, while in the second
experiment the cell is treated with the MT-depolymerizing agent nocodazole (NCZ) to
create a MT-less environment.
{\em In silico} we can easily turn on and off the advection term and compare cargo accumulation in the nucleus.
Notice that this corresponds to permitting or not the association to the MTs, while depolymerization 
cancels the MTs thus changing  as well the cell environment.
The full consequences on the cell metabolism under depolymerization are unknown,
and thus comparison with {\em in silico} experimentation is crucial.

\begin{figure}[h]
  \begin{center}
       \includegraphics[width=5.3cm]{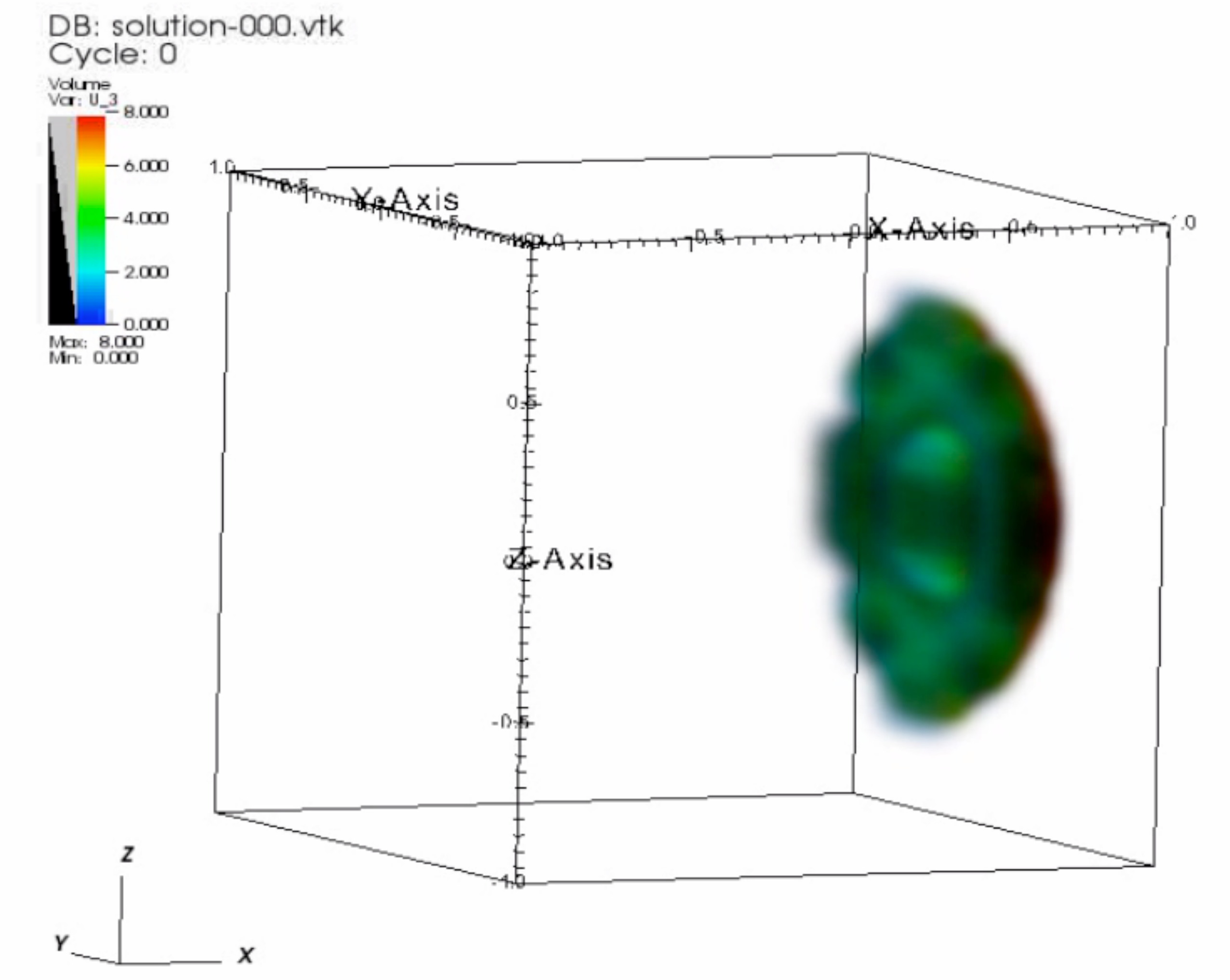}
\hspace{1 cm}
       \includegraphics[width=5.3cm]{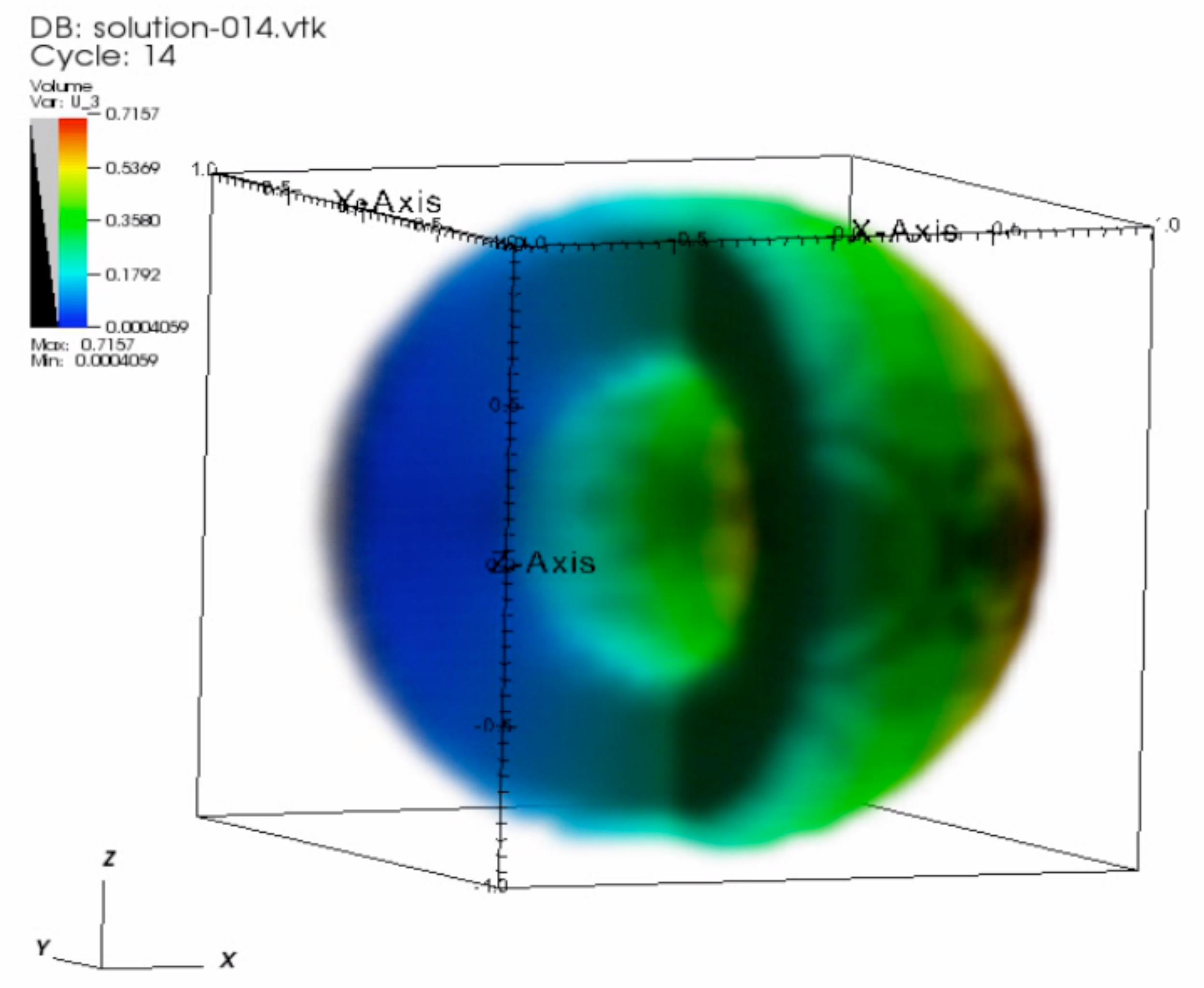}
  \end{center}
\caption {\footnotesize NLS cargo concentration in the cell 
at the initial time (left) and after 17 seconds (right).}
\label{fig:3d}
\end{figure}

The initial conditions are the concentrations of the single molecular species when cell is at rest (i.e. in the absence of external stimuli). 
Initial concentrations used in the simulation are reported in Table~\ref{tab:conc}. In particular we assume that the activated
cargo is initially concentrated in a peripheral zone of the cytoplasm,
as shown in the left plot in Figure~\ref{fig:3d}.
In such a situation, the concentrations of the complexes are zero (i.e. no complex is formed before the stimulus activates the cargo). The experimentally observed RanGTP 
accumulation of RanGTP in the nucleus forms in the initial phase of simulation. 

\begin{figure}[h]
  \begin{center}
       \includegraphics[width=9cm]{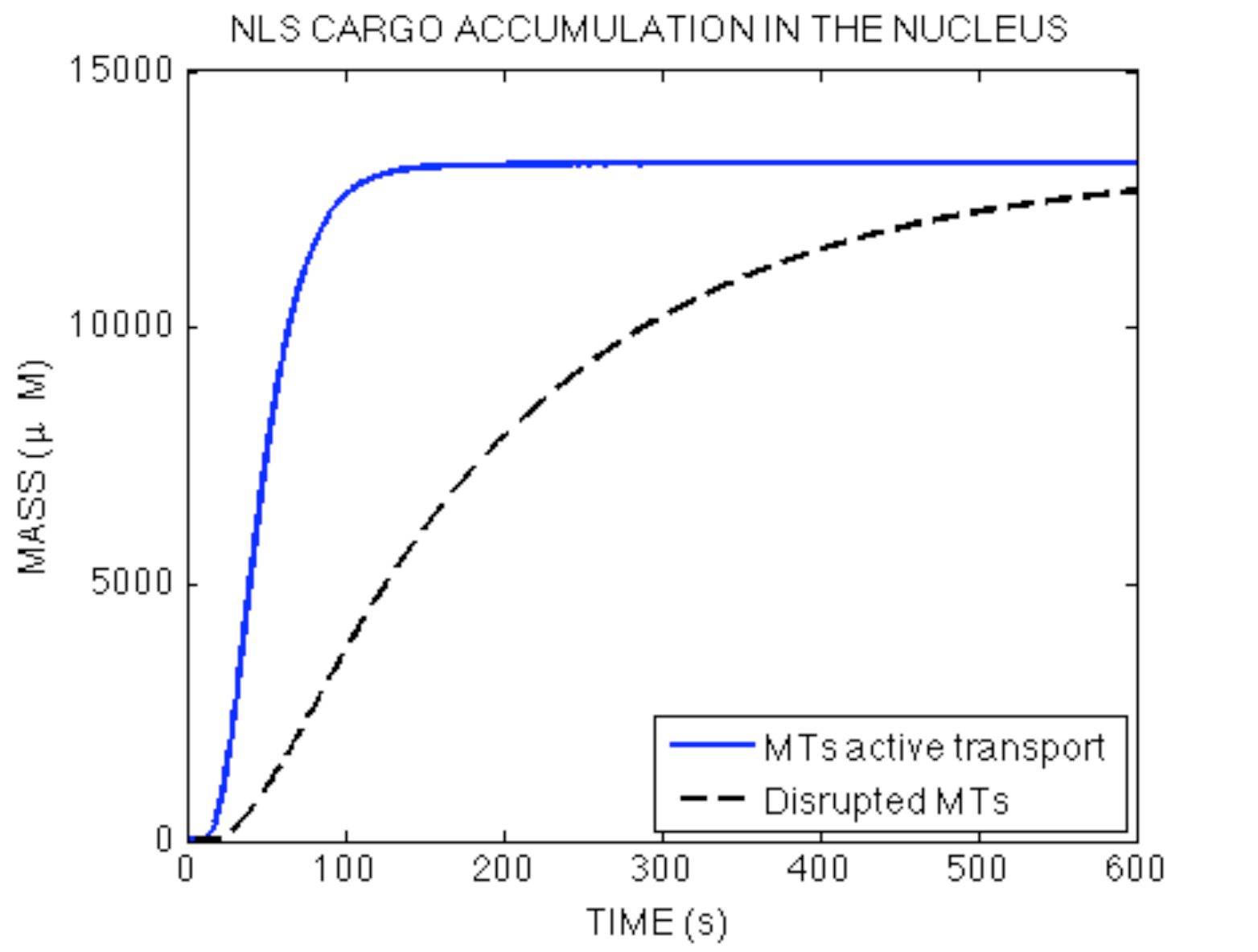}
  \end{center}
\caption {\footnotesize Accumulation of the cargo in the nucleus.}
\label{fig:cargoi}
\end{figure}

The accumulation in time of cargo mass in the nucleus is shown in Figure~\ref{fig:cargoi}. 
The change in accumulation rate is evident, and
confirms the behavior  experimentally obtained in~\cite{Roth07}.

\begin{figure}[h]
  \begin{center}
       \includegraphics[width=5.3cm]{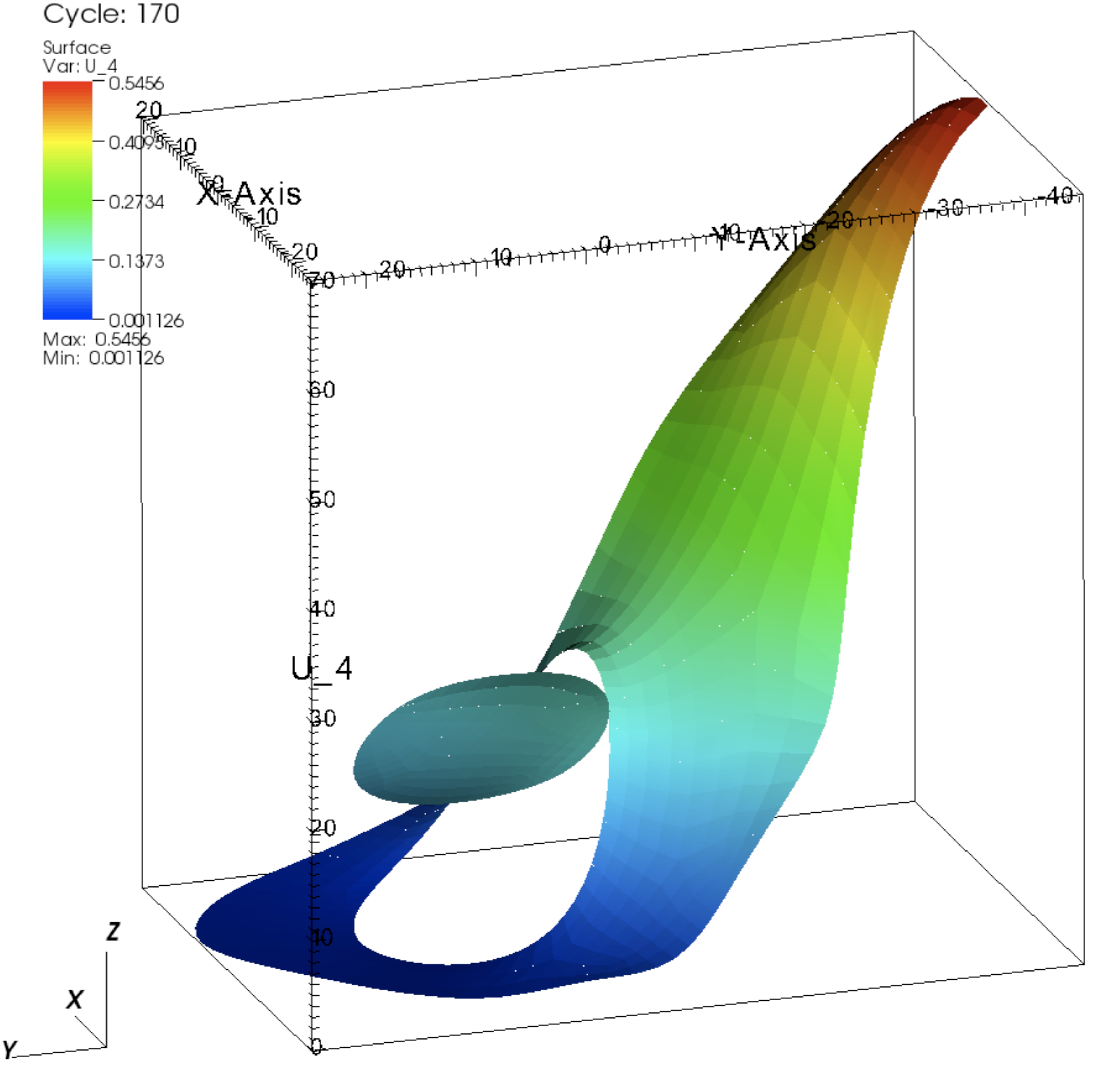}
\hspace{1 cm}
       \includegraphics[width=5.3cm]{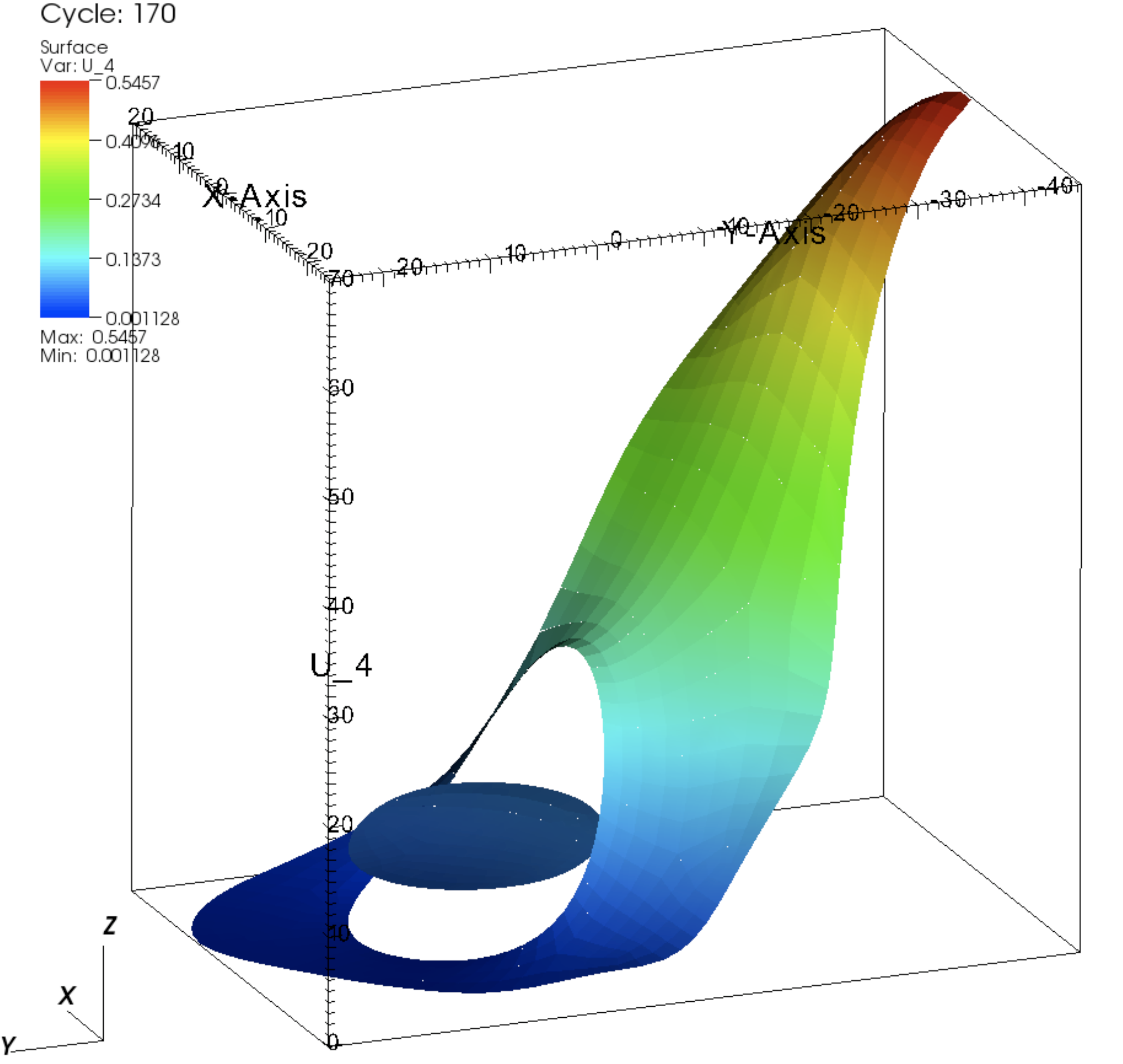}
  \end{center}
\caption {\footnotesize NLS cargo concentration in the cell 
after 17 second. With advection along the MTs(left) and without (right).}
\label{fig:surface}
\end{figure}

A snapshot showing the simulated cargo concentration in a spherical cell after 17 seconds is shown in Figure~\ref{fig:3d} (right), where
the accumulation jump across the NE is clearly visible. 
Figure~\ref{fig:surface} shows the concentration after the same time lapse, this time obtained with a two dimensional computation to better appreciate the variation of concentration along the cell. Notice in particular
that the concentration inside the nucleus is, in presence MTs active transport,  more
than twofold. 

\nocite{*}
\bibliographystyle{abbrv}
\addcontentsline{toc}{chapter}{Bibliography}
\bibliography{quaderno}

\begin{thebibliography}{10}

\bibitem{Alberts}
B.~Alberts, A.~Johnson, J.~Lewis, M.~Raff, K.~Roberts, and P.~Walter.
\newblock {\em Molecular Biology of the Cell (5th ed.)}.
\newblock Garland Science, 2008.

\bibitem{Axelrod76}
D.~Axelrod, D.~E. Koppel, J.~Schlessinger, E.~Elson, and W.~W. Webb.
\newblock Mobility measurement by analysis of fluorescence photobleaching
  recovery kinetics.
\newblock {\em Biophys J.}, 16(9):1055--1069, 1976.

\bibitem{becs05}
A.~Becskei and I.~Mattaj.
\newblock Quantitative models of nuclear transport.
\newblock {\em Curr. Opin. Cell Biol.}, 17(1):27--34, 2005.

\bibitem{Berez04}
A.~M. Berezhkovskii, Y.~A. Makhnovskii, M.~I. Monine, V.~Y. Zitserman, and
  S.~Y. Shvartsman.
\newblock Boundary homogenization for trapping by patchy surfaces.
\newblock {\em The Journal of Chemical Physics}, 121(22):11390--11394, 2004.

\bibitem{Berez06}
A.~M. Berezhkovskii, M.~I. Monine, C.~B. Muratov, and S.~Y. Shvartsman.
\newblock Homogenization of boundary conditions for surfaces with regular
  arrays of traps.
\newblock {\em The Journal of Chemical Physics}, 124(3):036103, 2006.

\bibitem{Briggs25}
G.~Briggs and J.~Haldane.
\newblock A note on the kinetics of enzyme action.
\newblock {\em Biochem. J.}, 197:338--339, 1925.

\bibitem{Camp03}
E.~M. Campbell and T.~J. Hope.
\newblock Role of the cytoskeleton in nuclear import.
\newblock {\em Advanced Drug Delivery Reviews}, 55(6):761--771, 2003.

\bibitem{Cangi09DG}
A.~Cangiani, E.~H. Georgoulis, and M.~Jensen.
\newblock Discontinuous galerkin methods for convection-diffusion-reaction
  systems with transmission conditions.
\newblock In preparation.

\bibitem{Cangi09m}
A.~Cangiani and R.~Natalini.
\newblock Models of cell signal transduction including active transport along
  the microtubules.
\newblock In preparation.

\bibitem{Caudron05}
M.~Caudron, G.~Bunt, P.~Bastiaens, and E.~Karsenti.
\newblock Spatial coordination of spindle assembly by chromosome-mediated
  signaling gradients.
\newblock {\em Science}, 309:1373--1376, 2005.

\bibitem{CockburnNASAVKI99}
B.~Cockburn.
\newblock Discontinuous {G}alerkin methods for convection-dominated problems.
\newblock In T.~Barth and H.~Deconink, editors, {\em High-{O}rder {M}ethods for
  {C}omputational {P}hysics}, volume~9, pages 69--224. Springer, 1999.

\bibitem{cks}
B.~Cockburn, G.~Karniadakis, and C.-W. Shu, editors.
\newblock {\em Discontinuous {G}alerkin Methods. {T}heory, Computation and
  Applications}, volume~11 of {\em DGMref}.
\newblock Springer, 2000.

\bibitem{deVries06}
G.~de~Vries, T.~Hillen, M.~Lewis, J.~Müller, and B.~Schönfisch.
\newblock {\em A Course in Mathematical Biology: Quantitative Modeling with
  Mathematical \& Computational Methods}.
\newblock SIAM, 2006.

\bibitem{Dinh06}
A.-T. Dinh, C.~Pangarkar, T.~Theofanous, and S.~Mitragotri.
\newblock Theory of spatial patterns of intracellular organelles.
\newblock {\em Biophysical Journal}, 90(10):67--69, 2006.

\bibitem{Dinh05}
A.-T. Dinh, T.~Theofanous, and S.~Mitragotri.
\newblock A model for intracellular trafficking of adenoviral vectors.
\newblock {\em Biophys. J.}, 89(3):1574--1588, 2005.

\bibitem{Eungd04}
N.~J. Eungdamrong and R.~Iyengar.
\newblock Modeling cell signaling networks.
\newblock {\em Biology of the Cell}, 96(5):355 -- 362, 2004.

\bibitem{Aebi03}
B.~Fahrenkrog and U.~Aebi.
\newblock The nuclear pore complex: nucleocytoplasmic transport and beyond.
\newblock {\em Nat Rev Mol Cell Biol}, 4(10):757--766, 2003.

\bibitem{Verkman91}
K.~Fushimi and A.~Verkman.
\newblock Low viscosity in the aqueous domain of cell cytoplasm measured by
  picosecond polarization microfluorimetry.
\newblock {\em J. Cell Biol.}, 112(4):719--725, 1991.

\bibitem{Gianna02}
P.~Giannakakou, M.~Nakano, K.~C. Nicolaou, A.~O'Brate, J.~Yu, M.~V.
  Blagosklonny, U.~F. Greber, and T.~Fojo.
\newblock Enhanced microtubule-dependent trafficking and p53 nuclear
  accumulation by suppression of microtubule dynamics.
\newblock {\em Proceedings of the National Academy of Sciences of the United
  States of America}, 99(16):10855--10860, 2002.

\bibitem{Gorlich03}
D.~G{\"o}rlich, M.~Seewald, and K.~Ribbeck.
\newblock Characterization of ran-driven cargo transport and the rangtpase
  system by kinetic measurements and computer simulation.
\newblock {\em EMBO J.}, 22:1088--1100, 2003.

\bibitem{Gunder99}
G.~G. Gundersen and T.~A. Cook.
\newblock Microtubules and signal transduction.
\newblock {\em Current Opinion in Cell Biology}, 11(1):81--94, 1999.

\bibitem{Harding05}
A.~Harding, T.~Tian, E.~Westbury, E.~Frische, and J.~Hancock.
\newblock Subcellular localization determines map kinase signal output.
\newblock {\em Curr Biol}, 15(9):869--873, 2005.

\bibitem{houston-etal:01}
P.~Houston, C.~Schwab, and E.~S{\"u}li.
\newblock Discontinuous {$hp$}-finite element methods for
  advection-diffusion-reaction problems.
\newblock {\em SIAM J. Numer. Anal.}, 39(6):2133--2163 (electronic), 2002.

\bibitem{ElbaumMT08}
A.~Kahana, G.~Kenan, M.~Feingold, M.~Elbaum, and R.~Granek.
\newblock Active transport on disordered microtubule networks: The generalized
  random velocity model.
\newblock {\em Physical Review E (Statistical, Nonlinear, and Soft Matter
  Physics)}, 78(5):051912, 2008.

\bibitem{Kao93}
H.~Kao, J.~Abney, and A.~Verkman.
\newblock {Determinants of the translational mobility of a small solute in cell
  cytoplasm}.
\newblock {\em J. Cell Biol.}, 120(1):175--184, 1993.

\bibitem{Khol06}
B.~Kholodenko.
\newblock Cell-signalling dynamics in time and space.
\newblock {\em Nat Rev Mol Cell Biol}, 7(3):165--176, 2006.

\bibitem{Klebe95a}
C.~Klebe, F.~Bischoff, H.~Ponstingl, and A.~Wittinghofer.
\newblock Interaction of the nuclear gtp-binding protein ran with its
  regulatory proteins rcc1 and rangap1.
\newblock {\em Biochemistry}, 34:639--647, 1995.

\bibitem{Klebe95b}
C.~Klebe, H.~Prinz, A.~Wittinghofer, and R.~Goody.
\newblock The kinetic mechanism of ran-nucleotide exchange catalyzed by rcc1.
\newblock {\em Biochemistry}, 34:12543--12552, 1995.

\bibitem{Elbaum07}
R.~B. Kopito and M.~Elbaum.
\newblock Reversibility in nucleocytoplasmic transport.
\newblock {\em Proceedings of the National Academy of Sciences},
  104(31):12743--12748, 2007.

\bibitem{Suli07}
A.~Lasis and E.~S\"{u}li.
\newblock \$hp\$-version discontinuous galerkin finite element method for
  semilinear parabolic problems.
\newblock {\em SIAM J. Numer. Anal.}, 45(4):1544--1569, 2007.

\bibitem{Maca01}
I.~G. Macara.
\newblock Transport into and out of the nucleus.
\newblock {\em Microbiol. Mol. Biol. Rev.}, 65(4):570--594, 2001.

\bibitem{Mattaj98}
I.~W. Mattaj and L.~Englmeier.
\newblock Nucleocytoplasmic transport: The soluble phase.
\newblock {\em Annual Review of Biochemistry}, 67(1):265--306, 1998.

\bibitem{MM13}
L.~Michaelis and M.~Menten.
\newblock Die kinetik der invertinwirkung.
\newblock {\em Biochem. Z.}, 49:333--369, 1913.

\bibitem{MurrayI}
J.~D. Murray.
\newblock {\em Mathematical Biology: I. An Introduction}.
\newblock Springer, 2002.

\bibitem{Nedel01}
F.~m.~c. N\'ed\'elec, T.~Surrey, and A.~C. Maggs.
\newblock Dynamic concentration of motors in microtubule arrays.
\newblock {\em Phys. Rev. Lett.}, 86(14):3192--3195, Apr 2001.

\bibitem{Parti98}
A.~Partikian, B.~Olveczky, R.~Swaminathan, Y.~Li, and A.~Verkman.
\newblock Rapid diffusion of green fluorescent protein in the mitochondrial
  matrix.
\newblock {\em J. Cell Biol.}, 140(4):821--829, 1998.

\bibitem{Quimby03}
B.~Quimby and M.~Dasso.
\newblock The small gtpase ran: interpreting the signs.
\newblock {\em Curr Opin Cell Biol}, 15(3):338--344, 2003.

\bibitem{Gorlich01}
K.~Ribbeck and D.~G{\"o}rlich.
\newblock Kinetic analysis of translocation through nuclear pore complexes.
\newblock {\em EMBO J.}, 20:1320--1330, 2001.

\bibitem{Riddick05}
G.~Riddick and I.~Macara.
\newblock A systems analysis of importin-$\alpha$-$\beta$ mediated nuclear
  protein import.
\newblock {\em J. Cell Biol.}, 168:1027--1038, 2005.

\bibitem{Riddick07}
G.~Riddick and I.~Macara.
\newblock The adapter importin-$\alpha$ provides flexible control of nuclear
  import at the expense of efficiency.
\newblock {\em Mol. Syst Biol.}, 3:118, 2007.

\bibitem{Roth07}
D.~M. Roth, G.~W. Moseley, D.~Glover, C.~W. Pouton, and D.~A. Jans.
\newblock A microtubule-facilitated nuclear import pathway for cancer
  regulatory proteins.
\newblock {\em Traffic}, 8:673--686(14), June 2007.

\bibitem{ElbaumMT05}
H.~Salman, A.~Abu-Arish, S.~Oliel, A.~Loyter, J.~Klafter, R.~Granek, and
  M.~Elbaum.
\newblock Nuclear localization signal peptides induce molecular delivery along
  microtubules.
\newblock {\em Biophysical Journal}, 89(3):2134--2145, 2005.

\bibitem{SS89}
L.~Segel and M.~Slemrod.
\newblock The quasi-steady-state assumption: A case study in perturbation.
\newblock {\em SIAM Review}, 31:446--477, 1989.

\bibitem{Seksek97}
O.~Seksek, J.~Biwersi, and A.~Verkman.
\newblock Translational diffusion of macromolecule-sized solutes in cytoplasm
  and nucleus.
\newblock {\em J. Cell Biol}, 138:1315--1342, 1997.

\bibitem{Smith02}
A.~Smith, B.~Slepchenko, J.~Schaff, L.~Loew, and I.~Macara.
\newblock Systems analysis of ran transport.
\newblock {\em Science}, 295:488--491, 2002.

\bibitem{Smith01}
D.~Smith and R.~Simmons.
\newblock Models of motor-assisted transport of intracellular particles.
\newblock {\em Biophys. J.}, 80(1):45--68, 2001.

\bibitem{Swamin96}
R.~Swaminathan, S.~Bicknese, N.~Periasamy, and A.~Verkman.
\newblock Cytoplasmic viscosity near the cell plasma membrane: translational
  diffusion of a small fluorescent solute measured by total internal
  reflection-fluorescence photobleaching recovery.
\newblock {\em Biophysical Journal}, 71(2):1140--1151, 1996.

\bibitem{Tekotte02}
H.~Tekotte and I.~Davis.
\newblock Intracellular mrna localization: motors move messages.
\newblock {\em Trends in Genetics}, 18(12):636 -- 642, 2002.

\end{thebibliography}
\end{document}